\documentclass[11pt]{article}

\usepackage{graphicx}
 
\newcommand{\pic}[2]{\raisebox{-.3\height}{\includegraphics[scale=#2]{#1}}}

 \usepackage{latexsym}
\newtheorem{Theorem}{Theorem}[section]
\newtheorem{Corollary}[Theorem]{Corollary}
\newtheorem{Proposition}[Theorem]{Proposition}
\newtheorem{Lemma}[Theorem]{Lemma}

\newcommand\BWtwoa{\pic{}{.700}}
\newcommand\BWtwob{\pic{}{.700}}
\newcommand\BWtwoc{\pic{}{.700}}
\newcommand\BWtwod{\pic{}{.700}}
\newcommand\BWtwoe{\pic{}{.700}}
\newcommand\BWtwof{\pic{}{.700}}
\newcommand\BWtwog{\pic{}{.700}}
\newcommand\BWtwoh{\pic{}{.700}}

\newcommand\BWfoura{\pic{}{.700}}
\newcommand\BWfourb{\pic{}{.700}}
\newcommand\BWfourc{\pic{}{.700}}
\newcommand\BWfourd{\pic{}{.700}}

\newcommand\BWfivea{\pic{}{.700}}
\newcommand\BWfiveb{\pic{}{.700}}
\newcommand\BWfivec{\pic{}{.700}}
\newcommand\BWfived{\pic{}{.700}}
\newcommand\BWfivee{\pic{}{.700}}
\newcommand\BWfivef{\pic{}{.700}}
\newcommand\BWfiveg{\pic{}{.700}}

\newcommand\BWsixa{\pic{}{.700}}

\newenvironment{Def}{\par\smallskip%
\noindent\textbf{Definition.}\  }%
{\par\smallskip}

\newenvironment{Rem}[1][{}]{\par\smallskip%
\noindent\textbf{Remark.}\ #1\ }%
{\par\smallskip}

\newenvironment{Notation}{\par\smallskip%
\noindent\textbf{Notation.}\  }%
{\par\smallskip}

\newenvironment{Proof}[1][{}]{\par\smallskip%
\noindent\textit{Proof #1: }\  }
{\hfill$\Box$\par\smallskip}

\newcommand{\be}{\begin{equation}}
\newcommand{\ee}{\end{equation}}

\begin{document}

\title{A basis for the Birman-Wenzl algebra}
\author{H. R. Morton\\[2mm]
University of Liverpool}
 
\date{}
\maketitle
\begin{abstract}
This paper   provides an explicit isomorphism between
the {\em Birman-Wenzl algebra} $BW_n$, constructed by J.\ Birman
and H.\ Wenzl, and the {\em Kauffman algebra}
$MT_n$, subsequently constructed by H.R.\ Morton and P.\
Traczyk. The Birman-Wenzl algebra is defined
algebraically using generators and relations while the Kauffman
algebra has a geometric formulation in terms of tangles. The
  isomorphism is obtained by constructing an explicit basis in
$BW_n$, analogous to a basis previously constructed for $MT_n$
using `Brauer connectors'. The geometric isotopy arguments used
for $MT_n$ are systematically replaced by algebraic versions using the
Birman-Wenzl relations. This not only gives a direct way of
determining the dimension of the Birman-Wenzl algebra, but also
clarifies the role played by the ring of coefficients,
$\Lambda$, and its specialisations. 
\end{abstract}

\section*{Foreword}

This paper is a very lightly edited version of an 
article originally written in 1989 but never fully completed.  It was intended to be a joint paper with A.J.Wassermann. He had planned to write
  a final section to make use of the ability to
change at will between the Birman-Wenzl algebra, as given by
generators and relations, and the geometric framework of the
tangles, so as to look in more detail at the representation
theory. 

One goal of our original approach was to make sure that
specialisations of the coefficient ring could be handled
confidently, and that the translations to and from the tangle
context were on a sound footing.

Subsequently others have made
progress in this way, in works such as
\cite{FG}, but we had a number of requests for our earlier
account, and so I put it into this more accessible form  on the Liverpool knot theory pages in 2000.  

In order to have a more permanent place for it I have now put it on ArXiv. I had hoped to make this   a joint submission, but I have been unable to re-establish contact with Wassermann to get his formal agreement. The present paper is largely the result of our joint discussions, although I take responsibility for the eventual content and exposition.

\section{Introduction}

In recent years there has been considerable interest in
deformations of the classical `centraliser algebras' of Schur,
Weyl and Brauer. These play an important role in several areas,
including exactly solvable models in statistical mechanics,
quantum groups, von Neumann algebras and knot theory. It has
long been recognized that these links are more than tenuous and
if properly exploited lead to fruitful interchanges between the
different disciplines. The first and most spectacular instance
of this was of course Vaughan Jones' pioneering work on
subfactors, which led to his discovery of new link invariants.
Subsequently these invariants were understood in terms of
solutions of the quantum Yang-Baxter equation and vertex models.
The central thread running through all these topics is the
quantum group obtained by deforming the universal enveloping
algebra of the unitary group. One has also to deform the
centraliser algebras, which amounts to replacing the group
algebra of the symmetric group by the Hecke algebra (of type
$A$). 

After these discoveries, somewhat curiously history took a
reverse turn. Kauffman discovered new link invariants, based on
a purely skein theoretic characterisation of Jones' original
invariants. It was natural to ask whether Kauffman's invariants
could be obtained by algebraic means. This led Birman and Wenzl
to introduce a deformation of an abstract algebra first
introduced by R.\ Brauer. An alternative knot-theoretic approach
to deforming Brauer's algebra was later given in terms of
tangles by Morton-Traczyk and by Kauffman himself. The original
algebra of Brauer bore the same relation to the orthogonal group
as the group algebra of the symmetric group did to the unitary
group. 

By exploiting the new insights provided by Vaughan Jones'
work on subfactors (in particular his `basic construction'),
Wenzl was able to acquire a fuller understanding of Brauer's
algebra, and resolve some old questions on semisimplicity raised
by Brauer and Weyl. Subsequent studies have shown that the
Birman-Wenzl algebra does indeed  provide the correct analogue
of the Hecke algebra for the quantum group corresponding to the
orthogonal group. Most recently Wenzl has been able to construct
new examples of subfactors using these algebras as a substitute
for the Hecke algebras.

In this paper we provide an explicit isomorphism between
the {\em Birman-Wenzl algebra} $BW_n$, constructed by J.\ Birman
and H.\ Wenzl in \cite{BW}, and the {\em Kauffman algebra}
$MT_n$, subsequently constructed by the author and P.\
Traczyk in \cite{MT}. The Birman-Wenzl algebra is defined
algebraically using generators and relations while the Kauffman
algebra has a geometric formulation in terms of tangles. We
obtain this isomorphism by constructing an explicit basis in
$BW_n$, analogous to a basis previously constructed for $MT_n$
using `Brauer connectors'. The geometric isotopy arguments used
in \cite{MT} are systematically replaced by algebraic versions using the
Birman-Wenzl relations. This not only gives a direct way of
determining the dimension of the Birman-Wenzl algebra, but also
clarifies the role played by the ring of coefficients,
$\Lambda$, an integral domain. 

In fact, in \cite{BW} the authors
prefer to consider the algebra $BW_n\otimes_{\Lambda}k$, where
$k$ is the field of fractions of $\Lambda$. This enables them to
imitate V.\ Jones' basic construction and thus determine the
structure of the algebra. At a crucial point in proof of their
main result (theorem 3.7) they need to use a specialisation of
$\Lambda$. Since the algebra is defined by generators and
relations over $\Lambda$, any such specialisation automatically
extends to $BW_n$ although not necessarily to the algebra
$BW_n\otimes_{\Lambda}k$. This difficulty can be overcome by
observing that the existence of a basis implies that $BW_n$ is
free as a module over $\Lambda$. The arguments presented in
\cite{Bourbaki} p.55  to prove that the Hecke algebra is generically
semisimple may then be adapted to prove the same result for
$BW_n$, i.e.\ the specialisations of $BW_n$ are semisimple for a
Zariski open subset of the parameter space
$\mbox{Spec}(\Lambda)$. Wenzl has carried out a more detailed
analysis, based on Jones' basic construction, in order to
determine precisely when the algebras fail to be semisimple.

This paper is divided into six sections. In section 2 we review
the definitions of the algebras to be studied, with some
historical comments. In section 3 we use the basic solution of
the Yang-Baxter equation for the orthogonal group together with
a simple skein-theoretic argument to provide a short
self-contained definition of Kauffman's two-variable link
invariant. We also briefly discuss the duality between the
quantum orthogonal group and the Birman-Wenzl algebra. In
section 4 we give an inductive definition of a basis for the
Birman-Wenzl algebra and outline the more formal aspects of the
proof. The inductive procedure relies on a natural filtration
analogous to the one extensively used by Hanlon and Wales in
their studies of Brauer's algebra \cite{HW}. Effectively the
proof that the proposed basis is a spanning set is achieved by a
double induction, which from the point of view of tangles
depends both on the number of strings and then on the number of
`through' strings. The remaining two sections are devoted to
various stages of the inductive argument showing that the
natural surjective maps from the Birman-Wenzl algebras to the
tangle algebras are isomorphisms. In section 5 we treat the case
in which there are no `through' strings: a complete
understanding of this case is crucial for the subsequent
reasoning  since it allows us to use geometry in place of
algebra in a controlled way. Finally in section 6 we perform the
main step of the induction.

\section{Three  algebras}

\subsection{The Birman-Wenzl algebra}
We start by recalling the definition of the Birman-Wenzl
algebra. We have made a slight change by the introduction of
some minus signs, in accordance with Kauffman's `Dubrovnik'
version of his link invariant. As explained in \cite{MT} and
below, this makes it much easier to see the connection with
Brauer's centraliser algebras.

Let $\Lambda$ be the quotient ring ${\bf
Z}[\lambda^{\pm1},z,\delta]/<\lambda^{-1}-\lambda-z(\delta-1)>$.
Thus
$\Lambda$ (or more accurately its complexification) is the
coordinate ring of the irreducible quasiprojective variety
defined by $\lambda\neq 0, \lambda^{-1}-\lambda=z(\delta-1)$ in
${\bf A}^3$. 
\begin{Def} The {\em Birman-Wenzl algebra} $BW_n$ is the
quotient of the free algebra over $\Lambda$ with generators
$g_1^{\pm1},g_2^{\pm1},\ldots,g_{n-1}^{\pm1}$ and
$e_1,e_2,\ldots,e_{n-1}$ modulo the ideal generated by the
relations:

\noindent(1)\ (Kauffman skein relation)\quad
$g_i-g_i^{-1}=z(1-e_i)$.

\noindent(2)\ (Idempotent relation)\quad
$e_i^2=\delta e_i$.

\noindent(3)\ (Braid relations)\quad
$g_ig_{i+1}g_i=g_{i+1}g_ig_{i+1}$ and $g_ig_j=g_jg_i$ if $\vert
i-j\vert >1$.

\noindent(4)\ (Tangle relations)\quad
  $e_i e_{i\pm1} e_i=e_i$ and $g_i g_{i\pm1} e_i =e_{i\pm1} e_i$.

\noindent(5)\ (Delooping relations) $g_ie_i=e_ig_i=\lambda e_i$
and
$e_ig_{i\pm1}e_i=\lambda^{- 1}e_i$.

\end{Def}

\begin{Rem} If $z$ is taken to be invertible then the idempotent
relation follows from the delooping and skein relations.

In Birman and Wenzl's original version several of their
relations could be omitted without loss, given invertibility of
$z$. They use $v$ in place of $\lambda $ in the coefficient ring.

The presentation given here is intended to be sufficiently
symmetric to allow for easy comparison with the tangle algebra,
while maintaining the coefficient ring $\Lambda$  as
in \cite{MT}.

\end{Rem}

\subsection{Kauffman's tangle algebra}

\begin {Def}  An $(m,n)\mbox {-tangle}$ is  a piece of
knot diagram in a rectangle $R$ in the  plane, consisting of
arcs and closed curves, so  that the end points of the arcs
consist of  $m$  points at the top of the rectangle and $n$
points  at the bottom, in some standard  position.\end{Def}  
An example of a  $(4,2)\mbox {-tangle}$ is shown in  figure 2.1.

\begin{center}
\BWtwoa\\
{Figure 2.1}
\end{center}

\medskip\nobreak 
\begin {Def}  Two tangles are {\em ambient  isotopic} if
they are related by a sequence of  Reidemeister's moves I, II
and III, (see figure  2.2), together with isotopies of $R$
fixing its  boundary.

\begin{center}
\BWtwob\\
{Figure 2.2}
\end{center}

\par \noindent  They are {\em regularly
isotopic} if Reidemeister  move I is not used.\end{Def} 
  
\begin{Notation}  Write  ${\cal U}_{n}^{m}$    for the
set of $(m,n)\mbox {-tangles}$ up to  regular
isotopy.\end{Notation} 

The set  ${\cal U}_{n}^{n}$ admits an 
{\em associative multiplication}, defined by  placing
representative tangles one below the  other. 

A well-known subset $B_{n}$ consists of geometric  braids, in
this context represented by tangles  (necessarily without closed
components) where the  height coordinate in $R$ increases
monotonically  on each component. It can be shown that $B_{n}$ 
is the full group of units in ${\cal  U}^{n}_{n}$ under the
multiplication.

The {\em closure},  $\hat T,$ of an 
$(n,n)\mbox {-tangle } T,$ is defined, by analogy  with the
closure of a braid, to be the link  diagram (or $(0,0)\mbox
{-tangle}$) given from $T$  by joining the points on the top of
$R$ to those  on the bottom by arcs lying outside $R$ with no 
further crossings.

We define  a closure map $\varepsilon :{\cal 
U}_{n}^{n}\to {\cal U}_{0}^{0}$, by 
$\varepsilon (T)=\hat T$.

\medskip  From  ${\cal U}_{n}^{n}$ we construct the 
algebra $MT_{n}$,   which we call Kauffman's tangle algebra, as
an algebra over a ring $\Lambda $, as in \cite{MT}.
We shall take $\Lambda $ to be the ring $$\Lambda  ={\bf
Z}[\lambda ^{\pm 1},z,\delta ]/<\lambda ^{-1}-\lambda =z(\delta
-1)>.$$
  Then $\Lambda $ is isomorphic to a subring of 
${\bf Z}[\lambda ^{\pm 1},z^{\pm 1}]$, by taking 
$\delta =1+(\lambda ^{-1}-\lambda )/z$. It admits  a
homomorphism $e:\Lambda \to {\bf Z}[\delta ]$ with
$e(z)=0,\  e(\lambda )=1$ and 
$e(\delta )=\delta .$ The main aim of this paper is to show that
the Birman-Wenzl algebra $BW_n$ is isomorphic to $MT_n$, on
specialisation of coefficients.
  
 Certain features of $MT_{n} $, for example its dimension, and 
its relation to Brauer's algebra \cite{B},  appear here very
simply, using the  homomorphism $e$ and the Dubrovnik invariant 
${\cal D}.$ These features of $BW_n$, not proved directly in the
original approach, then follow at once.
\medskip\nobreak 
\begin {Def}   {\em Kauffman's tangle algebra}, 
$MT_{n}$, is the $\Lambda $-module, constructed  from $\Lambda
[{\cal U}_{n}^{n}]$ by factoring  out three sets of relations:
\be T^{+}-T^{-}=z(T^{0}-T^{\infty }),\ee
where $T^{\pm },T^{0},T^{\infty }$ are  represented by tangles
differing only as in  figure 2.3,

\begin{center}
\BWtwoc\\
{Figure 2.3}
\end{center}

\be T^{\mbox {right}}=\lambda  ^{-1}T,\ 
T^{\mbox {left}}=\lambda T,\ee
where $T^{\mbox {right}}$ and $ T^{\mbox {left}} $  are given from
$T$ by adding a left or right hand  curl as in figure 2.4,
\begin{center}
\BWtwod\\ \nobreak
{Figure 2.4}
\end{center}

   \be T\amalg O=\delta T, \ee 
where $ T\amalg O$ is the union of $T$ and a  circle having no
crossings with $T$ or  itself.\end{Def}   

\begin {Proposition}  Composition of  tangles induces a
$\Lambda $-bilinear  multiplication on $MT_{n}$ making $MT_{n}$
an  algebra over $\Lambda$. \end{Proposition}  
\begin {Proof}  The relations (1)-(3) carry down  under the
multiplication in $\Lambda [{\cal  U}^{n}_{n}]$. \end{Proof} 
\begin {Proposition}  The map $\varepsilon 
$ induces a  $\Lambda $-linear map $\varepsilon  :MT_{n}\to
MT_{0}$. \end{Proposition}

We now give the homomorphism $\varphi  :BW_{n}\to MT_{n}$
which provided the  intuition behind Birman and Wenzl's
description  of $BW_{n}$. 
\begin {Def}  Write $G_{i}$,  $E_{i}$  respectively for
the tangles in ${\cal  U}_{n}^{n}$ illustrated in figure 2.5.
Use the  same letters for the elements represented by  these
tangles in $MT_{n}$,  called $s_{i}$,  
$h_{i}$ in \cite{MT}.\end{Def}

\begin{center}
\BWtwoe\\
{Figure 2.5}
\end{center}

   Then
$$G_{i}-G_{i}^{-1}=z(1-E_{i})$$ in $MT_{n}$,  from relation (1)
applied to the only crossing in 
$G_{i}$. 

Similarly, relation (2) shows that
\begin{eqnarray*}
	G_{i}E_{i}&=&E_{i}G_{i}=\lambda  ^{-1}E_{i} \\
   G_{i}^{-1}E_{i}&=&E_{i}G_{i}^{-1}=\lambda E_{i}
\end{eqnarray*} 
and relation (3) that  
\[
E_{i}^{2}=\delta
E_{i}.\]

\begin {Theorem}  A homomorphism $\varphi  :BW_{n}\to
MT_{n}$ may be defined by $\varphi  (g_{i})=G_{i},\ 
\varphi (e_{i})=E_{i}$. 
\end{Theorem}   
\begin {Proof}  The relations in $BW_{n}$ are  respected. We
have already noted that the skein  relation and delooping
relations are satisfied by 
$E_{i}$,  $G_{i}$ in $MT_{n}$. The other  relations hold even at
the level of the tangle  semigroup ${\cal U}_{n}^{n}$.
\end{Proof}

  Our goal is to prove that $\varphi $ is an 
isomorphism for all $n$. In this section we find  explicit
spanning sets for $MT_{n}$, and show  that  $\varphi $ is
surjective.

In section 4 we give the proof from \cite{MT}  that the chosen
spanning sets are a free basis  for $MT_{n}$, using the
existence of Kauffman's  invariant.

The proof that $\varphi $ is injective will  subsequently be
built up in stages, with the  recurring pattern of taking
spanning sets for  selected subspaces of $BW_{n}$ and proving
that  they map to independent sets in $MT_{n}$. 

To save later effort we note here some symmetry  of $BW_{n}$,
which carries over by $\varphi $ to  two natural operations in
$MT_{n}$. 
\begin {Def}  (1)\quad Write $\rho  _{n}:BW_{n}\to
BW_{n}$ for the automorphism  defined by $$\rho
_{n}(g_{i})=g_{n-i},\   \rho  _{n}(e_{i})=e_{n-i}.$$ 
\par \noindent (2)\quad Write $\alpha  :BW_{n}\to BW_{n}$
for the reversing {\em anti} automorphism defined by $$\alpha
(g_{i})=g_{i},\  
\alpha (e_{i})=e_{i}.$$ \end{Def}
\begin {Rem}  The symmetry of the relations  in $BW_{n}$
ensures that  $\rho _{n}$,  $\alpha $  are
well-defined.\end{Rem}
\begin {Proposition}  There is an automorphism $\rho
_{n}$ of $MT_{n}$, and an antiautomorphism $\alpha $, with
$\varphi \circ 
\alpha =\alpha \circ \varphi $ and $\varphi \circ 
\rho _{n}=\rho _{n}\circ \varphi $. \end{Proposition}   
\begin {Proof}  Write $\rho _{n},\, \alpha :{\cal
U}_{n}^{n}\to {\cal  U}_{n}^{n}$ for the natural
symmetries given by  rotating a tangle $T$ through $\pi $ about
one of the two axes shown in figure 2.6.

\begin{center}
\BWtwof\\
{Figure 2.6}
\end{center}

Clearly $\alpha (G_{i})=G_{i}$,  $\rho  _{n}(G_{i})=G_{n-i}$,
and similarly for $E_{i}$.  The skein relations are preserved by
$\rho _{n}$  and $\alpha $ so that they induce $\rho  _{n},\, \alpha :MT_{n}\to MT_{n}$.  Since $\rho _{n}(ST)=\rho
_{n}(S)\rho _{n}(T)$  and $\alpha (ST)=\alpha (T)\alpha (S)$
these are  respectively an automorphism and an 
antiautomorphism, satisfying the stated relations  on the
generators of $BW_{n}$. \end{Proof}
  We now continue with the
proof that $MT_{n}$ has  a finite spanning set, and at the same
time we  develop the notation to relate these algebras  readily
with Brauer's centraliser algebras.

\subsection {Connectors and Brauer's algebras}   An
$(n,n)$-tangle $T$ consists of $n$ arcs and a  number, $ \left
\vert T\right \vert  $, of closed  curves. If each arc joins a
point at the top to a  point at the bottom then the tangle
determines a  permutation in $S_{n}$.

\begin{Def}
For a general tangle we extend the idea of a  permutation to
that of an $n$-{\em connector}, defined  to be a pairing of $2n$
points into $n$ pairs.
\end{Def}
 The set $C_{n}$ of $n$-connectors has 
$(2n)!/2^{n}n!$ elements, the product of the  first $n$ odd
integers.

Take the set of $2n$ points to be the end points  of  $(n,n)\mbox
{-tangles}.$ The arcs of any  $T\in  {\cal U}_{n}^{n}$ pair
these end points to give  a connector, which we write as $\mbox
{conn}(T)\in  C_{n}.$

\begin{Rem}  [(Brauer's algebra)]  Brauer \cite{B} uses 
$C_{n}$  as the basis for an algebra over  ${\bf  Z}[\delta
],$  (writing  $ n $   in place of   
$\delta $   and   $f$   in place of   $n$ ).      He divides
the  $2n$   points to be connected   into two subsets  
$t_{1},\dots{} ,t_{n}$   and   
$b_{1},\dots{} ,b_{n},$ arranged along the top  and bottom of a
rectangle, and views a connector  
$c$ as a set of  $n$   intervals with these  $2n$  points as
endpoints, which join the points paired  by $c.$      Two
connectors   $c_{1}$   and   
$c_{2}$   are composed by placing one rectangle  above the
other, giving   $n$   arcs whose  endpoints are the new top and
bottom points,  together with some number   $r\geq 0$   of
closed  curves. 

   Brauer sets   $c_{1}c_{2}=\delta ^{r}d,$    where   $d$   is
the connector defined by the new  arcs.    This defines an
associative  multiplication on   ${\bf Z}[\delta 
][C_{n}]=A_{n}$   making it an algebra over   
${\bf Z}[\delta ],$ called {\em Brauer's algebra.  }

Having divided the $2n$ points in this way there  is a natural
embedding $S_{n}\subset C_{n}.$ 
\end{Rem}   

We can modify the map $\mbox {conn}:{\cal  U}_{n}^{n}\to
C_{n}$ to give a   multiplicative homomorphism $c:{\cal 
U}_{n}^{n}\to A_{n},$ which extends to 
$c:MT_{n}\to A_{n}$ as follows. 

For $T\in {\cal U}_{n}^{n}$ set $c(T)=\delta ^{ 
\left \vert T\right \vert }\mbox {conn}(T)\in  A_{n}$. This can be
extended to $c:\Lambda  [{\cal U}_{n}^{n}]\to A_{n}$ by
setting 
$c(\Sigma \, \lambda _{i}T_{i})=\Sigma 
\, e(\lambda _{i})c(T_{i})$, using the  ring
homomorphism $e:\Lambda \to {\bf  Z}[\delta ]$. 
\begin {Theorem}  There is an induced  homomorphism
$c:MT_{n}\to A_{n}.$ \end{Theorem}
\begin {Proof}  The relations (1)-(3) defining 
$MT_{n}$ are respected.\end{Proof} 
\begin {Rem}  We show later that $A_{n}$ is  isomorphic to
the  ${\bf Z}[\delta ]$ algebra 
$MT_{n}\otimes _{\Lambda }{\bf Z}[\delta ]$  given from $MT_{n}
$ by replacing the  coefficients $\Lambda $ with ${\bf
Z}[\delta ],$  using the homomorphism $e.$ 

The existence of $c:MT_{n}\to A_{n}$ can be  viewed as the
consequence of specialising the  coefficients so that the
relations no longer  distinguish under- from over-crossings.
Then  tangles pass to their projections, retaining only  the
information of their connectors. The crucial  technical feature
here is that we can specialise 
$\Lambda $ so as to retain $\delta ,$ while  fixing $\lambda =1 $
and $z=0$. Complications arise if  we try to do this while
working in the ring
${\bf Z}[\lambda ^{\pm  1},z^{\pm 1}].$ \end{Rem}

\begin {Def}  Given a tangle $T,$ choose  a sequence of
base-points, consisting firstly of  one end point of each arc,
and then one point on  each closed component. Say that $T$ is 
{\em totally descending } (with this choice of  base points) if
on traversing all the strands of 
$T,$ starting from the base point of each  component in order,
each crossing is first met as  an overcrossing.\end{Def}
  
\begin {Rem}   We shall assume that for each connector $c\in 
C_{n}$ a choice of  ordering of base-points for  the arcs has
been made, and we use this same  choice for all tangles $T$ 
with  $c=\mbox  {conn}\, T$. Note that there are  
$n!\, 2^{n}$ potentially different  choices possible for
each connector. The precise  choice is not material, and we
shall have  occasion to vary the choice in the course  of
later proofs. The result will be simply to  alter the choice of
linear basis in $MT_{n}$.  
\end{Rem}

An example of a totally descending $(3,3)\mbox  {-tangle}$ is
shown in figure 2.7, with  base-points numbered according to a
choice of  order.

\begin{center}
\BWtwog\\
{Figure 2.7}
\end{center}

\begin {Theorem}\label{2.6}  $MT_{n}$ is spanned by  totally
descending tangles.\end{Theorem}   
\begin {Proof}  Let $T$ be a tangle representing  an element of
$MT_{n}$. Choose base points for 
$T$ according to the choice for $\mbox {conn}(T)$.  Traverse the
arcs of  $T$ in order. At the first  non-descending crossing use
relation (1) with 
$T=T_{\pm }$. Note that $\mbox {conn}(T_{+})=\mbox  {conn}(T_{-})$,
so that $T_{\mp}  $, resulting  from $T$ with the crossing
switched, has fewer  non-descending crossings. Then  $T$ is a
linear  combination of three tangles, two with fewer  crossings
and one with fewer non-descending  crossings. The theorem
follows by induction,  firstly on the number of crossings, then
on the  number of non-descending crossings.\end{Proof} 
\begin {Corollary}   $MT_{n}$ is spanned by  totally
descending tangles without closed  components.\end{Corollary}  

\begin {Proof}  If $T$ is totally descending,  with $r$ closed
components, then these components  are unknotted curves stacked
below the arcs of 
$T$. The tangle can then be altered by regular  isotopy so that
the unknotted components lie well  away from the arcs.  By using
Reidemeister move I  as well they can be changed to have no 
self-crossings. Then by (2) and (3),  $T=\lambda  ^{k}\delta
^{r}T'$ in $MT_{n},$ where $T'$  consists simply of the arcs of
$T.$ \end{Proof} 
\begin {Rem}  This result holds as stated for 
$n=0$, provided that we admit the `empty tangle'  as an element
of ${\cal U}^{0}_{0}$. In any  event $MT_{0}$ is spanned by a
single  element.\end{Rem}   
\begin {Corollary}   $MT_{0}$ is  cyclic.\end{Corollary}  

\begin {Theorem} \label{2.9} Let $S$ and $T$ be  totally descending
$(n,n)\mbox {-tangles},$ without  closed components, such that
$\mbox {conn}(S)=\mbox  {conn}(T).$ Then $S$ and $T$ are ambient 
isotopic, and so $S=\lambda ^{k}T$ in  $MT_{n}$,  for some $k$.
\end{Theorem}  
\begin {Proof}  Number the arcs of  $S$ and $T$  according to
the order of their base points.  Since $\mbox {conn}(S)=\mbox
{conn}(T),$ the $i\mbox  {th}$ arc in each tangle joins the same
pair of  end points. The arcs can be arranged to lie in 
disjoint levels $1$ to $n$ above the plane of 
$R,$ since arc $i$  lies above arc $j$ at every  crossing when
$i<j.$ Each individual arc is  unknotted, because the tangle is
descending, so  it can be changed by ambiemt isotopy to an arc 
without self-crossings in its level. The  resulting tangles are
then ambient isotopic by  level-preserving isotopy.\end{Proof}  
\begin {Rem}  If the arcs of  $S$ and $T$  have no
self-crossings initially then $S$ and $T$  are regularly
isotopic.\end{Rem}   

\begin{Rem} [(Construction)]  For each connector $c\in  C_{n},$
choose an order for the arcs. With this  order construct a
totally descending tangle with  connector $c$ such that any two
arcs cross at  most once. (Start for example from a diagram of 
the connector in which any two arcs cross at most  once, and
make it descending, by choosing the  sense of each crossing
according to the order of  the arcs.) The element $T_{c}\in
{\cal  U}_{n}^{n}$ represented by this tangle then  depends only
on $c$ and the chosen order, by  Theorem \ref{2.9}.\end{Rem}   
\begin {Rem}  For $c\in S_{n}$  and a natural  choice of
order the resulting tangles $T_{c}$  have been studied,
\cite{E,EM}, under the name  `positive permutation braids'. They
can be  represented by a braid in $B_{n}$ with positive 
crossings and permutation $c$ in which any two  strings cross at
most once.

These braids have also been used in \cite{MS1,MS2},  to give
easily handled generators for the Hecke  algebra $H_{n}.$
\end{Rem} 

\begin {Theorem}   $MT_{n}$ is spanned, for  every choice
of order, by the finite set $\lbrace  T_{c}\rbrace ,\  c\in
C_{n}$. \end{Theorem}   
\begin {Proof}  By theorem \ref{2.6} and its corollary, 
$MT_{n}$ is spanned by tangles which are ambient  isotopic to 
$T_{c},$ for various $c.$ By use of  relation (2), any tangle
ambient isotopic to 
$T_{c}$ represents $\lambda ^{k}T_{c}$ in 
$M_{n},$ for some $k.$ \end{Proof} 
\begin {Rem}  The number of crossings in a  totally
descending tangle $T_{c}$ depends on the  connector $c$, not on
the order of arcs used. It  is simply the number of pairs of
arcs which cross  in $c$, as dictated by whether or not their 
endpoints interlock on the boundary rectangle.

Clearly any tangle with $k$ crossings can always  be written in
$MT_{n}$  as a linear combination  of totally descending tangles
with at most $k$  crossings, by induction on $k$, using the 
procedure of theorem \ref{2.6}. It follows that if 
$T'_{c}$,  $T_{c}$ are totally descending tangles  with the same
connector $c$, arising from  different choices of the order of
arcs then
$$
	T'_{c}=T_{c}+\sum _{d} \lambda _{d}T_{d},
$$ where $d$ runs over connectors with fewer  crossings than
$c$. \end{Rem}

  We finish this section by proving:
\begin {Theorem}  The map $\varphi  :BW_{n}\to MT_{n}$
is surjective.\end{Theorem}

\begin {Proof}  We must show that $MT_{n}$ is  generated by
$E_{i}$,  $G_{i},\  1\leq i\leq  n-1$. It is enough to show
that each totally  descending tangle $T_{c}$ is a monomial in 
$\lbrace E_{i}\rbrace $ and $\lbrace G_{i}^{\pm  1}\rbrace $. 

Assuming that the connector $c$ pairs $r$ points  at the top
with $r$ at the bottom, and connects  the remaining $2k=n-r$
points as $k$ pairs,  leaving $2k$ points at the bottom
connected as  
$k$ pairs.

We can then draw the tangle $T_{c}$ (for any  order of the arcs)
so that there are $r$ arcs  running monotonically from top to
bottom, $k$  arcs running with a single local minimum from top 
to top, and $k$ arcs from bottom to bottom with a  single local
maximum. We can further assume,  since the arcs never cross
twice, that all the  local minima on the top arcs are higher up
than  the local maxima, so that there are only  $r$  arcs
passing through the middle part of the  rectangle.

Now pair arbitrarily the local maxima and minima,  and isotop
the tangle so that each local minimum  moves down to lie
directly above its  corresponding maximum. We can now decompose
the  tangle level by level into a composite of simple  tangles
in each of which there are $n$ strings  all running vertically,
except for one pair,  which either cross simply, giving
$G_{i}^{\pm  1}$, or form a paired minimum and maximum, giving 
$E_{i}$. An example is shown in figure 2.8.

\begin{center}
\BWtwoh\\
{Figure 2.8}
\end{center}

\end{Proof} 
\begin {Rem}  It is useful to regard the  tangle $T_{c}$
with $r$ through strings as a  composite of an $(n,r)$-tangle
and an 
$(r,n)$-tangle, and it suggests that a  counterpart of
$(n,r)$-tangles might helpfully be  studied in relation to
$BW_{n}$. \end{Rem}

\section {Kauffman's link polynomial}
In this section we discuss Kauffman's Dubrovnik 
invariant of links, and its relation to the solutions of 
the Yang-Baxter equation for the orthogonal 
group.

Kauffman's polynomial, in its Dubrovnik form, is 
a non-zero function ${\cal D}:{\cal 
U}_{0}^{0}\to \Lambda ,$   i.e. a function on 
knot diagrams which is unaltered by {\em regular } 
isotopy.

\medskip This function ${\cal D}$ has 
the basic properties:
\par \noindent (1)\quad $
{\cal D}(K^{+})\  -\  {\cal 
D}(K^{-})\  \  =\  \  
z({\cal D}(K^{0})\  -\  {\cal 
D}(K^{\infty }))\qquad \mbox {(skein relation)}$ 

 where the diagrams $K^{\pm },\  K^{0}$ and 
$K^{\infty }$ differ only as in figure 2.3, and
 \par \noindent  (2)\quad$ {\cal D}(K^{\mbox 
{left}})=\lambda {\cal D}(K),\  \  
{\cal D}(K^{\mbox {right}})=\lambda ^{-1}{\cal 
D}(K),$ 

 where $K^{\mbox {left}}$ and $K^{\mbox {right}}$ are 
given from $K$ as in figure 2.4.

It also satisfies
 \par \noindent (3)\quad $ {\cal D}(K\amalg 
O)=\delta {\cal D}(K), $

 where $K\amalg O$ is the union of $K$
and  a circle having no crossings with $K$ or with 
itself, and $\delta \in \Lambda $ satisfies 
$\lambda ^{-1}-\lambda =z(\delta -1).$

\begin {Proposition}  Kauffman's invariant 
exists if and only if the cyclic module  $MT_{o}$ 
is free.\end{Proposition}   
\begin {Proof}  We have shown already that 
$MT_{0}$ is cyclic, so $MT_{0}$ is free if and 
only if there is a non-zero  $\Lambda 
$-homomorphism $\varphi :MT_{0}\to \Lambda $.

If $MT_{0}$ is free then we may define ${\cal 
D}$ on any diagram $K$ by ${\cal D}(K)=\varphi 
(K)$.
Conversely, if ${\cal D}$ satisfies (1)-(3) 
then it defines a non-zero $\Lambda 
$-homomorphism ${\cal D}:MT_{0}\to \Lambda 
$.\end{Proof} 
\begin{Rem} [(Uniqueness of Kauffman's invariant)]  
It follows simply from section 2 that Kauffman's 
invariant is unique, because $MT_{0}$ is cyclic. 
It is determined uniquely by its value on $O$, 
the diagram of the unknot without any crossings. 
 ${\cal D}$  was originally normalised so that 
${\cal D}(O)=1$.  It now appears more natural 
to assign the value 1 to the `empty knot', so 
that ${\cal D}(O)=\delta $.\end{Rem}   

Kauffman's original proof of the existence of 
${\cal D}$, \cite{K}, requires a considerable 
amount of combinatorial argument to show that the 
elements of $\Lambda $ reached by different 
routes from a given diagram $K$ are independent 
of any intermediate choices.

We note here an alternative existence proof, 
using the Yang-Baxter orthogonal invariants.

\begin {Proposition} There exists a regular 
isotopy invariant of knot diagrams in  ${\bf 
Z}[s^{\pm 1}]$ which satisfies relations (1)-(3) 
with $z=s-s^{-1},\  \lambda 
=s^{2n-1},\  \delta =1+(\lambda -\lambda 
^{-1})/z$, and takes the value 1 on the empty 
knot.\end{Proposition}   
\begin {Proof}[\rm{(Turaev)}] The invariant is 
constructed from the $q$-analogue of the 
fundamental representation of the Lie algebra of 
$SO(2n)$. \end{Proof} 

For each $n$ we have a ring homomorphism 
$e_{n}:\Lambda \to {\bf Z}[s^{\pm 1}]$ 
defined by $e_{n}(\lambda )=s^{2n-1},\  
e_{n}(z)=s-s^{-1}$. Turaev's invariant then 
defines a map $\varphi _{n}:MT_{0} \to   
{\bf Z}[s^{\pm 1}]$ with $\varphi 
_{n}(aK)=e_{n}(a)\varphi _{n}(K)$ for $a\in 
\Lambda $.

\begin {Proposition}   $MT_{0}$ is a free 
$\Lambda $-module.\end{Proposition}  
\begin {Proof}  Suppose not. Then there exists 
$a\in \Lambda ,\  a\mathbin{\not =}0$ such 
that $aK=0$, where $K$ is the empty diagram.
Now $\varphi _{n}(K)=1$ so $0=\varphi 
_{n}(aK)=e_{n}(a)$ for all $n$. This is 
impossible, since for any given  $a\mathbin{\not 
=}0$ there exists $n$ with $e_{n}(a)\mathbin{\not 
=}0$. \end{Proof} 

This proves the existence of ${\cal D}$, given 
Turaev's invariants.  In principle ${\cal 
D}(K)$  could be calculated explicitly for a 
given link diagram $K$ from knowledge of the invariants 
$\varphi _{n}(K)$ for sufficiently many $n$, as 
follows:
\begin{Proof}
We know that any element  $a$ of $\Lambda $ can 
be written as a polynomial in $\lambda ^{\pm 1}$, 
$z$ and $\delta $. Now $z\delta =\lambda 
^{-1}-\lambda +z$, so   $z^{k}a$ can be rewritten 
as a polynomial in $\lambda ^{\pm 1}$ and $z$ 
alone, for large enough $k$.

A simple induction, as in theorem \ref{2.6}, shows that 
 $z^{ \left
\vert K\right \vert }{\cal  D}(K)\in \Lambda $ can always be
written as a polynomial in
$\lambda  ^{\pm 1}$ and $z$; say

\begin{eqnarray*} z^{ \left \vert K\right \vert 
}{\cal D}(K)&=&\sum _{r=m}^{M} \lambda 
^{r}P_{r}(z) \\
			&=&\sum _{r=m}^{M}\lambda ^{r}Q_{r}(s), 
\end{eqnarray*}
where $Q_{r}(s)=P_{r}(s-s^{-1})$.

It is then enough to find $Q_{r}(s),\,m\leq 
r\leq M$.

Now for each $n$,
\begin{eqnarray*}\sum 
_{r=m}^{M}s^{r(2n-1)}Q_{r}(s)&=&e_{n}( z^{ \left 
\vert K\right \vert }{\cal D}(K)) \\ 
				&=&(s-s^{-1})^{ \left \vert K\right \vert 
}\varphi _{n}(K).\end{eqnarray*}

Write $V$ for the  $k\times k$ Vandermonde matrix 
with entries $$s^{(2n-1)r},\  1\leq n\leq 
k,\  m\leq r\leq M,\quad \mbox {with 
}k=M-m+1.$$
Then 
$$V\pmatrix{Q_{m} \cr  Q_{m+1} 
\cr  \vdots   \cr  Q_{M}}= 
(s-s^{-1})^{ \left \vert K\right \vert 
}\pmatrix{\varphi _{1} \cr  \varphi 
_{2}  \cr  \vdots  \cr  \varphi 
_{k}}.$$
Since $V$ is invertible, we have $Q_{m},\dots{} 
,Q_{M}$, and hence ${\cal D}(K)$ in terms of 
$\varphi _{1},\dots{} ,\varphi _{k}$.\end{Proof} 

In order to make these calculations explicitly we 
need bounds for $m$ and $M$, in terms of  $K$. It 
is certainly sufficient to note that $ \left 
\vert m\right \vert ,M\leq  \left \vert K\right 
\vert +c(K)$, where  $c(K)$ is the number of 
crossings in the diagram, although these bounds 
may turn out to  be quite generous.

\section {A basis for the tangle algebra}

In this section we set out the induction to be 
used in proving that the algebra $BW_{n}$ defined 
by generators and relations is isomorphic to the 
Kauffman algebra, defined by tangles. We start by 
reviewing the position for  $MT_{n}$. 

The algebra $MT_{n}$ is shown in \cite{MT} to be 
free over $\Lambda $, of the same dimension, $ 
\left \vert C_{n}\right \vert $, as  Brauer's 
algebra $A_{n}$. The proof, which we give here, 
is an easy consequence of the existence, however 
established, of Kauffman's Dubrovnik invariant 
${\cal D}:MT_{0}\to \Lambda $.  We make use 
of the homomorphism $e:\Lambda \to {\bf 
Z}[\delta ]$.
\begin {Proposition} \label{4.1}  $e({\cal 
D}(K))=\delta ^{ \left \vert K\right \vert }$. 
\end{Proposition}   
\begin {Proof} It follows from condition (1) that 
$ e({\cal D}(K))$ is unaltered when any 
crossing in a diagram is switched, and from (2) 
that it is unaltered by Reidemeister move I. Now 
any  diagram can be changed to any other with the 
same number of components by a sequence of 
crossing switches and Reidemeister moves, so 
$e({\cal D}(K))= e({\cal D}(K'))$, where $K'$ 
is the disjoint union of $ \left \vert K\right 
\vert $ simple closed curves, giving the result 
by (3).  \end{Proof}  

\begin {Theorem}\label{4.2}  Any set of tangles 
$\lbrace T_{c}\rbrace ,\  c\in C_{n}$, 
without closed components, having $c=\mbox {conn}\, (T_{c})$ and spanning $MT_{n}$ 
forms a free $\Lambda $-basis for 
$MT_{n}$.\end{Theorem}   
\begin {Proof}  Define a bilinear map 
$b:MT_{n}\times MT_{n}\to \Lambda $ by 
$b(S,T)={\cal D}(\varepsilon (ST)).$ Write $A$ 
for the $ \left \vert C_{n}\right \vert \times  
\left \vert C_{n}\right \vert $ matrix  with 
entries $a_{cd}=b(T_{c},T_{d}).$ 

Suppose that $\Sigma \lambda _{i}T_{i} =0,\, \lambda _{i}\in \Lambda $. We want to show 
that $\lambda _{i}=0$ for all $i$. For each $c\in 
C_{n}$ replace the $c$th column of $A$ 
by the  linear combination of the columns of $A$ 
with coefficients $\lambda _{i}$.
The new matrix then has determinant $\lambda 
_{c}\mbox {det}\, A$ and a zero column.
The required result follows by proving that $\mbox {det}\, A\not =0$, since $\Lambda $ has 
no zero-divisors.

Now $\varepsilon (T_{c}T_{d})\in MT_{0}$ is 
represented by a link with $r$ components, say. 
Each component contains at least one arc from  
$T_{c}$ and one from $T_{d},$ so $r\leq n.$ When 
$r=n$ each component must have exactly one arc 
from each, so that the connector $d$ is the 
`mirror image' of $c,$ given by interchanging the 
roles of the top and bottom points. Set 
$\overline c=d$ in this case, so that we have 
$r=n$ if and only if $d=\overline c.$ 

Now apply the homomorphism $e$ to the entries in 
$A.$ Then, by proposition \ref{4.1},  $e(a_{cd})=\delta 
^{r},\  r\leq 
n,$ and $r=n$ if and only if $d=\overline c.$ 
The matrix $e(A)$ has then one entry $\delta 
^{n}$ in each row and column, so $e(\mbox {det}\, A)=\mbox {det}(e(A))\in \bf 
Z[\delta ]$ has a non-zero coefficient for 
$\delta ^{n^2}.$ Thus $e(\mbox {det}\, A)
\not =0,$ so  $\mbox {det}\, A\not =0.$ 
\end{Proof} 

This  shows that $MT_{n}$ is a 
deformation of Brauer's algebra $A_{n}$, in the 
following sense. 

\begin {Theorem}  There is an isomorphism 
of $\bf Z[\delta ]$-algebras induced by $c$ 
between $MT_{n}\otimes _{\Lambda }\bf Z[\delta 
]$ and $A_{n}.$ \end{Theorem}   
\begin {Proof}  The map $c:MT_{n}\to A_{n}$, 
defined in section 2, 
factors through a $\bf Z[\delta ]$-homomorphism 
$MT_{n}\otimes _{\Lambda }\bf Z[\delta ]\to 
A_{n}.$ Since $MT_{n}\otimes _{\Lambda }\bf 
Z[\delta ]$ is spanned over $\bf Z[\delta ]$ 
by $\lbrace T_{c}\rbrace $ which maps onto a 
{\em basis } of $A_{n}$ of the same cardinality, 
this set must be a $\bf Z[\delta ]$-basis in the 
specialisation, and the map is hence an 
isomorphism.\end{Proof} 
\begin {Corollary} [to theorem \ref{4.2}] \label{4.2a}
 Any set
of  tangles with distinct connectors forms an 
independent set in $MT_{n}$.\end{Corollary}   
\begin {Proof}  We have shown that $c: 
MT_{n}\to A_{n}$ carries a free $\Lambda 
$-basis to a free ${\bf Z}[\delta ]$-basis. It 
follows, using determinantal criteria for 
independence as in theorem \ref{4.2}, that $k$ elements of 
$MT_{n}$whose images are independent in $A_{n}$ 
must themselves be independent.\end{Proof} 

We shall prove by induction on $n$ that the 
homomorphism $\varphi :BW_{n}\to MT_{n}$ is 
an isomorphism. In the course of the proof we 
shall construct explicit bases $\varphi 
^{-1}\lbrace T_{c}\rbrace $ in  $BW_{n}$. As part 
of the induction we shall use  natural 
filtrations $BW_{n}^{(r)}$ and $MT_{n}^{(r)}$ by 
2-sided ideals, analogous to the filtration of 
$A_{n}$ used by Hanlon and Wales, \cite{HW}. In 
the case of $MT_{n}$ this filtration arises from 
the geometric viewpoint, as in \cite{MT}, when 
we consider tangles of rank $\leq r$.
\begin {Def}  A tangle $T\in {\cal 
U}_{n}^{n} $ has {\em rank } $\leq r$ if it is the 
composite  $T=AB$ of an $(n,r)$ and an $(r,n)$ 
tangle.\end{Def}   
\begin {Rem}  Then $\mbox {conn}(T)$ has at 
most $r$ arcs connecting top to bottom. However 
this is not sufficient for $T$ to have rank $r$. 
For example, the tangle $T$ in figure 4.1 has 
rank 2, although $\mbox {conn}(T)$ has no 
connecting arcs 
from top to bottom.\end{Rem}   

\begin{center}
\BWfoura\\
{Figure 4.1}
\end{center}

Write $MT_{n}^{(r)}$ for the subspace of $MT_{n}$ 
spanned by tangles of rank $\leq r$. Clearly 
$MT_{n}^{(r)}$ is a 2-sided ideal, with
$$
	MT_{n}=MT_{n}^{(n)}\supset  
MT_{n}^{(n-2)}\supset \dots{} \  .
$$

\begin {Proposition}   $MT_{n}^{(r)}$ is 
generated, as an ideal, by the element 
$E_{1}E_{3}{\dots}E_{2k-1}$, where 
$2k=n-r$.\end{Proposition}   
\begin {Proof}  For $r>0$ we can write the 
identity tangle in ${\cal U}_{r}^{r}$ as
$$I=C(E_{1}E_{3}{\dots}E_{2k-1})D,$$
where $C$ is an $(r,n)$ tangle and $D$ is an 
$(n,r)$ tangle, as in figure 4.2. 

\begin{center}
\BWfourb\\
{Figure 4.2}
\end{center}

Then any tangle  $T=AB$ of rank $\leq r$ can be 
written as  $T=AC(E_{1}E_{3}{\dots}E_{2k-1})BD$ 
with $AC,BD\in {\cal U}_{n}^{n}$.

The case $r=0$ can be handled similarly, by first 
writing a tangle $T$ of rank 0 as $T=AE_{1}B$ 
where $A$ is an $(n,2)$ tangle and $B$ is a 
$(2,n)$ tangle. \end{Proof} 
\begin {Def}  For $r=n-2k$ write 
$BW_{n}^{(r)}$ for the 2-sided ideal of $BW_{n}$ 
generated by $e_{1}e_{3}{\dots}e_{2k-1}$.%
\end{Def}   
Then
$$BW_{n}=BW_{n}^{(n)}\supset 
BW_{n}^{(n-2)}\supset \dots{} \  .$$
Clearly $\varphi :BW_{n}\to MT_{n}$ restricts 
to $\varphi :BW_{n}^{(r)}\to MT_{n}^{(r)}$.

Our main result, that $\varphi $ is an 
isomorphism, follows from
\begin {Theorem}\label{4.5}   $\varphi 
:BW_{n}^{(r)}\to MT_{n}^{(r)}$ is injective 
for all $n,r$.\end{Theorem}   
\begin {Proof}  The detailed lemmas needed appear 
in later sections. The scheme of the proof 
follows here.

For fixed $n$ we prove the result for $r=0,1$ in 
section 5 from the 
injectivity of $\varphi $ on $BW_{n-1}$   
($=BW_{n-1}^{(n-1)}$) using induction on $n$.

The proof then continues by induction on $r$.

 For this induction we construct a linear 
subspace $V_{n}^{(r)}\subset BW_{n}^{(r)}$, 
complementing $BW_{n}^{(r-2)}$. The induction 
step for injectivity of $\varphi $ follows by 
establishing:
\par \noindent (1)\quad    
$V_{n}^{(r)}+BW_{n}^{(r-2)}$ is a 2-sided ideal 
in $BW_{n}^{(r)},$
\par \noindent (2)\quad  $\varphi \vert 
V_{n}^{(r)}\to MT_{n}$ is injective,
\par \noindent (3)\quad  
$e_{1}e_{3}{\dots}e_{2k-1}\in V_{n}^{(r)}$.

In the construction, given later in this section, 
we exhibit an explicit spanning set for 
$V_{n}^{(r)}$ whose image in $MT_{n}$ is an 
independent set of totally descending tangles. 
This establishes property (2).

Property (3) is immediate from the construction, 
and property (1) is proved in section 6. 
\end{Proof}

To describe certain elements in $BW_{n}$ we now 
draw on Artin's braid group.

The braid group on $n$ strings, defined by 
geometric braids, (particular types of  $(n,n)$ 
tangles), is known to have the presentation with 
generators $\sigma _{i},\  i\leq n$ and 
relations $$\sigma _{i}\sigma _{j}=\sigma 
_{j}\sigma _{i},\   \left \vert i-j\right 
\vert >1,\qquad \sigma _{i}\sigma _{i+1}\sigma 
_{i}=\sigma _{i+1}\sigma _{i}\sigma _{i+1}.$$
There is then a homomorphism $\psi :B_{n}\to 
BW_{n}$ defined by $\sigma _{i}\mapsto g_{i}$.
 Any two monomials in $BW_{n}$ in $g_{i}^{\pm 1}$ 
which arise from the same geometric braid  $\beta 
$ will then be equal, and we shall use  $\beta $ 
to picture the element $\psi (\beta )$. We shall 
also refer to monomials in $g^{\pm 1}$  as braids 
in $BW_{n}$.

There is an antihomomorphism $\mbox {perm}:B_{n}\to S_{n}$ defined by $\mbox {perm}(\sigma _{i})=\tau _{i}=(i\  i+1)$. 
With our convention of composition of geometric 
braids, the strings in a braid $\beta $ then join 
the point $i$ at the top to the point $\pi (i)$ 
at the bottom, with $\pi =\mbox {perm}(\beta )$.

Among the elements of  $B_{n}$ we shall use 
particularly the {\em positive permutation braids 
} and, as special cases, the {\em Lorenz braids }.
\begin {Def}  A braid in $B_{n}$ in which 
all crossings are positive and every pair of 
strings crosses at most once is called a 
{\em positive permutation braid }.\end{Def}   
\begin {Theorem}  A positive permutation 
braid $\beta $ is determined by the permutation 
$\pi =\mbox {perm}(\beta )$ induced by its 
strings.\end{Theorem}   
\begin {Proof}  Such braids are examples of 
`totally descending tangles', as defined in 
section 2, in which the arcs of the connector all 
join top to bottom and are ordered by the order 
of their initial points.\end{Proof} 
We shall write $\beta _{\pi }$ for the positive 
permutation braid with permutation $\pi =\mbox {perm}(\beta _{\pi })$, whose strings join the 
points $i$ at the top with $\pi (i)$ at the 
bottom. The element $b_{\pi }=\psi(\beta _{\pi 
})\in BW_{n}$, which we shall also call a 
positive permutation braid, can be conveniently 
referred to by the permutation $\pi $, rather 
than choosing one 
of the many ways of writing it as a monomial in 
$g_{i}$. For example, the permutation $\pi 
=(14)(23)\in S_{4}$ gives $b_{\pi 
}=g_{1}g_{2}g_{3}g_{1}g_{%
2}g_{1}=g_{2}g_{1}g_{2}g_{3}g_{2}g_{1}={\dots}$.

\begin {Def}  A {\em Lorenz braid } of 
type $(\ell ,r)$ is a braid $\beta _{\pi }$ where 
$\pi \in S_{n},\  n=\ell +r$, does not 
permute the first $\ell $ `left-hand' strings, or 
the last $r$ `right-hand' strings among 
themselves.\end{Def}   
For fixed $(\ell ,r)$ there are $(^{n}_{r})$ 
Lorenz braids, as a Lorenz permutation $\pi $ is 
determined simply by the free choice of endpoints 
for the right-hand strings.   Note that $\pi $ is 
an $(\ell ,r)$ Lorenz permutation if and only if 
$\pi (i)<\pi (j)$ for $1\leq i<j\leq \ell $ and 
for $\ell +1\leq i<j\leq n$. An example of a 
$(3,4)$ Lorenz braid is shown in figure 4.3.

\begin{center}
\BWfourc\\
{Figure 4.3}
\end{center}

Where $\pi^{-1} $ is a Lorenz permutation the 
braid 
$\beta_{\pi}=\alpha(\beta _{\pi^{-1} })$ can be 
viewed as a Lorenz 
braid $\beta _{\pi^{-1} }$  turned upside down. 
Call  $\alpha (\beta _{\pi ^{-1}})$ a {\em reverse 
Lorenz braid }. Note 
that  $(\beta _{\pi })^{-1}$ is not the same 
braid as  $\beta _{\pi ^{-1}}$ but has all the 
crossings switched.

\begin {Def}  For each $n$ and $r=n-2k$ 
write $V_{n}^{(r)}$ for the linear subspace of 
$BW_{n}^{(r)}$ spanned by elements $b_{\pi 
}w_{2k}b_{\tau }b_{\mu }$, where $\pi ,\mu $ 
are $(2k,r)$ Lorenz permutations,  $\tau $ is a 
permutation of the last $r$ strings only, and 
$w_{2k}\in BW_{2k}^{(0)}$.\end{Def}   

\begin {Proposition}  Given that $\varphi 
\vert BW_{2k}^{(0)}$ is injective for $n\geq 
2k$ then $\varphi \vert V_{n}^{(r)}\to 
MT_{n}$ is injective.\end{Proposition}   
\begin {Proof}  We know that $MT_{2k}^{(0)}$ is 
spanned by $ \left \vert 
C_{k}\right \vert ^{2} $ totally descending 
tangles, one for each  $k$-connector of rank 0, 
and that $\varphi \vert BW_{2k}^{(0)}\to 
MT_{2k}^{(0)}$ is surjective. By hypothesis we 
can choose a spanning set of $ \left \vert 
C_{k}\right \vert ^{2} $ elements for 
$BW_{2k}^{(0)}$ with this set of tangles as 
image.

Then $V_{n}^{(r)}$ is spanned by the $ 
(^{n}_{r})^{2} \left \vert 
C_{k}\right \vert ^{2}r!$ elements $b_{\pi 
}w_{2k}b_{\tau }b_{\mu }$, where $\pi ^{-1},\, \mu $ are drawn independently from $(2k,r)$ 
Lorenz permutations, $\tau $ from permutations in 
$S_{r}$ and $w_{2k}$ from the spanning set for 
$BW_{2k}^{(0)}$. 

The elements $\varphi (b_{\pi }w_{2k}b_{\tau 
}b_{\mu })$ are represented by tangles in 
$MT_{n}^{(r)}$ each with exactly $r$ through 
strings, and all having different connectors. A 
typical such tangle with $k=2,\, r=4$ is 
illustrated in figure 4.4. It follows by the 
corollary to theorem \ref{4.2} that these tangles are 
independent in $MT_{n}$, and hence that $\varphi 
\vert V_{n}^{(r)}$ is injective. \end{Proof}

\begin{center}
\BWfourd\\
{Figure 4.4}
\end{center}

This establishes property (2) of theorem \ref{4.5} 
under its induction hypothesis.

From theorem \ref{4.5} we eventually build a basis for 
$BW_{n}$  as a union of 
spanning sets for each $V_{n}^{(r)}$. The image 
of 
this basis in $MT_{n}$ can be represented by a 
set of  tangles each with a different connector, 
and each totally descending, for some ordering of 
the arcs. 

We note that this gives a complicated check that 
the dimension of $BW_{n}$ is
$$ \left \vert C_{n}\right \vert =\sum 
_{k=0}^{[n/2]} (^{n}_{r})^{2} \left \vert 
C_{k}\right \vert ^{2}r!,$$ where we write 
$r=n-2k$.

\section {Generators and relations for the tangle algebra: the
base for induction}
 
In this section we prove injectivity of $\varphi 
$ on $BW_{n}^{(0)}$ or $BW_{n}^{(1)}$, depending 
on the parity of $n$, given injectivity of 
$\varphi \vert BW_{n-1}$. The corresponding 
sets of tangles in $MT_{n}$ are those with at 
most one through string.

We start with some results in $BW_{n}$ which use 
only the regular isotopy relations.
\begin {Def}  The shift map 
$S:MT_{n}\to MT_{n+1}$  is a homomorphism 
defined on an $n$-tangle  $T$ as shown in figure 
5.1.\end{Def}   

\begin{center}
\BWfivea\\
{Figure 5.1}
\end{center}

Thus $S(G_{i})=G_{i+1},\  S(E_{i})=E_{i+1}$.

It is clear, from the behaviour on tangles, as 
shown in figure 5.2, that $WA_{m}=A_{m}S(W)$ for 
$W\in MT_{m}$, where $A_{m}=G_{m}G_{m-1}\dots{} 
G_{1}$.

\begin{center}
\BWfiveb\\
{Figure 5.2}
\end{center}

We can define a shift map with similar properties 
in  $BW$ as follows.
\begin {Def}  The shift map 
$S:BW_{n}\to BW_{n+1}$ is defined as a 
homomorphism by $$S(g_{i})=g_{i+1},\  
S(e_{i})=e_{i+1},$$
extended linearly.\end{Def}   
It is simply necessary to check that the 
relations are respected by $S$.
\begin {Proposition}\label{5.1}  The homomorphism $S$ 
satisfies
$$wa_{m}=a_{m}S(w),\quad wb_{m}=b_{m}S(w)$$ for 
any $w\in BW_{m}$, where
$$a_{m}=g_{m}g_{m-1}\dots{} g_{1},\  
b_{m}=g_{m}^{-1}g_{m-1}^{-1}\dots{} 
g_{1}^{-1}.$$\end{Proposition}   
\begin {Proof}  When $w=g_{i}^{\pm 1}$ or 
$w=e_{i}$ the result is an immediate consequence 
of the relations, and it follows   for monomials 
$w$ by induction on their length.\end{Proof}

We now define $F_{k}\in MT_{n},\  2k\leq n$, 
to be the element represented by the tangle shown 
in figure 5.3.

\begin{center}
\BWfivec\\
{Figure 5.3}
\end{center}

The following equations in $MT_{n}$ are clear 
from inspection of representative tangles.
\begin {Proposition} \label{5.2} For all $i<k$ 
\par \noindent (1) \quad 
$G_{i}F_{k}=G_{2k-i}F_{k} $,
\par \noindent (2) \quad 
$E_{i}F_{k}=E_{2k-i}F_{k} $,
\par \noindent (3) \quad 
$F_{k}G_{i}=F_{k}G_{2k-i}$,
\par \noindent (4) \quad 
$F_{k}E_{i}=F_{k}E_{2k-i}$. \end{Proposition}  
An example of equation (1) is illustrated in 
figure 5.4, with $i=1$ and $k=3$.

\begin{center}
\BWfived\\
{Figure 5.4}
\end{center}

Again it is clear from inspection of the tangles, 
as shown in figure 5.5, that
$$F_{k}=\alpha (A_{2k-2})
F_{k-1}E_{2k-1}A_{2k-2}.$$
By analogy we define $f_{k}\in BW_{n}$ 
inductively, setting  $f_{0}=\mbox {identity}$, and
$$f_{k}=\alpha (a_{2k-2})
f_{k-1}e_{2k-1}a_{2k-2}.$$
We then have  $F_{k}=\varphi (f_{k})$, and 
$\alpha (f_{k})=f_{k}$, since $f_{k-1}$ and 
$e_{2k-1}$ commute. 

\begin{center}
\BWfivee\\
{Figure 5.5}
\end{center}

\begin {Rem}  While it is clear that $\rho 
_{2k}(F_{k})=F_{k}$ in $MT_{n}$, it is difficult 
to prove directly from the definition and 
relations in $BW_{n}$ that $\rho 
_{2k}(f_{k})=f_{k}$. \end{Rem}   
\begin {Proposition}   
$BW_{n}^{(n-2k)}\subset BW_{n}$ is the 2-sided 
ideal generated by  $f_{k}$.\end{Proposition}   
\begin {Proof}  By definition $BW_{n}^{(n-2k)}$ 
is the 2-sided ideal generated by 
$e_{1}e_{3}\dots{} e_{2k-1}$. By induction on $k$ 
we can write $f_{k}=\alpha 
(b_{k})e_{1}e_{3}\dots{} e_{2k-1}b_{k}$ for some 
invertible element $b_{k}\in BW_{2k}$, in fact  
$b_{k}$ can be chosen to be a braid.\end{Proof} 

We now make use of the relations in $BW_{n}$ to 
prove the analogous results to proposition \ref{5.2}.
\begin {Proposition}\label{5.4}  For all  $i<k$ 
\par \noindent (1)\quad  
$g_{i}f_{k}=g_{2k-i}f_{k}$,
\par \noindent (2)\quad  
$e_{i}f_{k}=e_{2k-i}f_{k}$,
\par \noindent (3)\quad 
$f_{k}g_{i}=f_{k}g_{2k-i}$,
\par \noindent (4)\quad 
$f_{k}e_{i}=f_{k}e_{2k-i}$.
\par \noindent The same results hold with $\rho 
_{2k}(f_{k})$ in place of $f_{k}$.\end{Proposition}   
\begin {Proof}  Cases (3) and (4) follow from (1) 
and (2) by applying  $ \alpha $.
Applying $\rho _{2k}$ gives the results for $\rho 
_{2k}(f_{k})$.
The result is immediate for $k=1$.
 For $i>1$ the result follows from proposition 
\ref{5.1} by induction on $k$.

For example, in case (1),
\begin{eqnarray*}g_{i}f_{k}&=&g_{i}\alpha 
(a_{2k-2})f_{k-1}e_{2k-1}a_{2k-2} \\  
&=&\alpha (a_{2k-2})g_{i-1
}f_{k-1}e_{2k-1}a_{2k-2},\mbox { by applying 
}\alpha \mbox { to \ref{5.1}} \\
&     =&\alpha (a_{2k-2})g_{2k-i-1}f
_{k-1}e_{2k-1}a_{2k-2},\mbox { by 
induction} \\ 
& =&g_{2k-i}\alpha (a_{2k-2})f_{k-
1}e_{2k-1}a_{2k-2},\  (i\geq 2) 
\\   
&=&g_{2k-i}f_{k}. 
\end{eqnarray*}
To prove \ref{5.4} when i=1 we  set  $h_{j}=\alpha 
(a_{j})\alpha (a_{j-2})e_{j+1}e_{j-1}$.

Since $f_{k}=h_{2k-2}f_{k-2}a_{2k-4}a_{2k-2}$, 
the result for cases (1) and (2) will follow by 
showing that
\par \noindent $(1')$\quad  
$g_{1}h_{j}=g_{j+1}h_{j}$
\par \noindent $(2')$\quad 
$e_{1}h_{j}=e_{j+1}h_{j}$,
\par \noindent for all $j$.

We prove $(1')$ and $(2')$ by induction on $j$, 
starting with $j=2$. For  $j=2$ we have
$$h_{2}=g_{1}g_{2}e_{1}e_
{3}=e_{2}e_{1}e_{3}=e_{2}
e_{3}e_{1}=g_{3}g_{2}e_{3}e_{1}.$$
Then $g_{3}h_{2}=g_{3}g_{
1}g_{2}e_{1}e_{3}=g_{1}h_{2}$ and 
$e_{1}h_{2}=e_{1}e_{2}e_{
1}e_{3}=e_{1}e_{3}=e_{3}h_{2}$.

For the induction step, use the braid relations 
to write
$$\alpha (a_{j})\alpha 
(a_{j-2})=g_{2}g_{1}S(\alpha (a_{j-1})\alpha 
(a_{j-3})).$$
(Compare the two braids illustrated in figure 
5.6.)

\begin{center}
\BWfivef\\
{Figure 5.6}
\end{center}

Then $h_{j}=g_{2}g_{1}S(h_{j-1})$. So
\begin{eqnarray*}g_{1}h_{j}&=&g_
{1}g_{2}g_{1}S(h_{j-1})=g
_{2}g_{1}g_{2}S(h_{j-1})=
g_{2}g_{1}S(g_{1}h_{j-1}) \\  
&=&g_{2}g_{1}S(g_{j}h_{j-1}),\mbox { by induction on 
}j, \\  
&=&g_{2}g_{1}g_{j+1}S(h_{j-1})=g_{j+1}h_{j},\mbox { 
for }j>2. 
\end{eqnarray*}
Similarly $e_{1}h_{j}=e_{j+1}h_{j}$, using the 
relation in $BW_{n}$  that 
$e_{1}g_{2}g_{1}=g_{2}g_{1}e_{2}$.\end{Proof} 
\begin {Lemma} \label{5.5} Suppose that $\varphi 
:BW_{m+1}\to MT_{m+1}$ is injective. Then 
$BW_{m+1}e_{m}=BW_{m}e_{m}$.\end{Lemma}   
\begin {Proof}  By hypothesis it is enough to 
prove the corresponding result 
$$MT_{m+1}E_{m}=MT_{m}E_{m} .$$

For an  $(m+1,m+1)$ tangle $T$ define  
$\varepsilon _{m}(T)$ to be the $(m,m)$ tangle 
shown in figure 5.7.

\begin{center}
\BWfiveg\\
{Figure 5.7}
\end{center}

Using the standard interpretation of $\varepsilon 
_{m}(T)$ as an $(m+1,m+1)$ tangle it is clear 
that $\varepsilon _{m}(T)E_{m}=TE_{m}$. Extend 
the definition of $\varepsilon _{m}$ to linear 
combinations of tangles to define a linear map 
$\varepsilon _{m}:MT_{m+1}\to MT_{m}$, (the 
relations are respected).
Then any element $XE_{m}$ with $X\in MT_{m+1}$ 
can be rewritten as $XE_{m}=\varepsilon 
_{m}(X)E_{m} \in MT_{m}E_{m}$.\end{Proof} 
\begin {Corollary}  Under the same conditions,  
$BW_{m+1}e_{1}=S(BW_{m})e_{1}$.\end{Corollary}   
\begin {Proof}  Apply the automorphism $\rho 
_{m+1}$.\end{Proof} 
\begin {Proposition}\label{5.6}  Suppose that $\varphi 
\vert BW_{n-1}\to MT_{n-1}$ is injective.
Then \par \noindent 
(1)\quad$BW_{2k}f_{k}=BW_{k}f_{k}$, for all $k$ 
with $2k\leq n$, 
\par \noindent (2)\quad 
$BW_{2k+1}S(f_{k})=BW_{k+1}S(f_{k})$, for all $k$ 
with $2k+1\leq n$. \end{Proposition}   
\begin {Proof} 

\noindent(1)\quad The case $k=1$ is 
immediate, since $g_{1}^{\pm 1}e_{1} $ and 
$e_{1}e_{1}$ are multiples of $e_{1}$.

For $k\geq 2$ we have $k+1\leq n-1$ so that 
$BW_{k+1}e_{k}=BW_{k}e_{k}$ by lemma \ref{5.5}. It is 
enough to show that
\begin{eqnarray*} g_{i}BW_{k}f_{k}&\subset &BW_{k}f_{k} 
\\  e_{i}BW_{k}f_{k}&\subset &
BW_{k}f_{k}, \mbox { for all }i<2k. 
\end{eqnarray*}
This is immediate for  $i<k$. For $i>k$ it 
follows from \ref{5.4}, and the fact that $BW_{k}$ then 
commutes with $e_{i}$ and $g_{i}$.

Write $f_{k}=e_{k}r_{k}$ for some  $r_{k}\in 
BW_{2k}$, by induction on $k$. The remaining 
cases with $i=k$ follow by noting that 
 \begin{eqnarray*}g_{k}BW_{k}e_{k}&\subset &
BW_{k+1}e_{k}=BW_{k}e_{k} \\  \mbox { 
and }e_{k}BW_{k}e_{k}&\subset &
BW_{k+1}e_{k}=BW_{k}e_{k}.
\end{eqnarray*} 
\par \noindent (2)\quad The case  $k=1$ will be 
proved directly.

For $k\geq 2$ we have $k+2\leq n-1$ so that 
$BW_{k+2}e_{k+1}=BW_{k+1}e_{k+1}$ from lemma \ref{5.5}.
We must show that
\begin{eqnarray*}g_{i}BW_{k+1}S(f_{k})&\subset &
BW_{k+1}S(f_{k}), \\ 
e_{i}BW_{k+1}S(f_{k})&\subset &
BW_{k+1}S(f_{k}),\mbox { for }i\leq 2k. 
\end{eqnarray*}
This is immediate for $i<k+1$. For $i>k+1$ it 
follows from proposition \ref{5.4}, since $BW_{k+1}$ 
commutes with $g_{i}$ and $e_{i}$. The remaining 
cases follow as in (1), since 
$S(f_{k})=e_{k+1}S(r_{k})$.

We finish the proof of (2) by showing that 
$BW_{3}e_{2}=BW_{2}e_{2}$. Now  $BW_{2}e_{2}$ is 
spanned by $e_{2},\, e_{1}e_{2}$ and 
$g_{1}e_{2}$, so we must show that products of 
these elements with  $g_{2}$ or $e_{2}$ on the 
left still lie in $BW_{2}e_{2}$.  It is a matter 
of a quick check from the relations in  $BW_{3}$, 
to see that $e_{2}^{2}=\delta e_{2}$,  
$e_{2}e_{1}e_{2}=e_{2}$, $e_{2}g_{1}e_{2}=\lambda 
^{-1}e_{2}$, $g_{2}e_{2}=\lambda e_{2} $, 
$g_{2}e_{1}e_{2}=g_{1}^{-1}e_{2}$ and  
$g_{2}g_{1}e_{2}=e_{1}e_{2}$.\end{Proof} 

\begin {Corollary}\label{5.6a}  Suppose that $\varphi 
\vert BW_{n-1}$ is injective. Then the ideals 
generated by $f_{k}$ in $BW_{n}$, with $k=[n/2]$, 
 can be written as: 
\par \noindent (1)\quad 
$BW_{2k}^{(0)}=BW_{k}f_{k}BW_{k}$ when $n=2k$, 
and
\par \noindent (2)\quad 
$BW_{2k+1}^{(1)}=BW_{k+1}S(f_{k})BW_{k+1}$ when  
$n=2k+1$.\end{Corollary} 
\begin {Proof}[(1)] 
\begin{eqnarray*} 
BW_{2k}^{(0)}&=&BW_{2k}f_{k}BW_{2k} \\  
&=&BW_{k}f_{k}BW_{2k}\mbox { by \ref{5.6},} \\  
&=&BW_{k}f_{k}BW_{k}\mbox { applying }\alpha \mbox { 
to \ref{5.6}.} 
\end{eqnarray*}
\end{Proof}
\begin{Proof}[(2)]
The ideal  
$BW_{2k+1}^{(1)}$ generated by $f_{k}$ is equally 
generated by $S(f_{k})=a_{2k}^{-1}f_{k}a_{2k}$ so 
the result follows using \ref{5.6} (2) exactly as in 
(1).\end{Proof} 

We complete this section by showing the 
injectivity of $\varphi $ on the 2-sided ideals 
generated by $f_{k}$ in $BW_{n},\  k=[n/2]$, 
given injectivity on  $BW_{n-1}$.
\begin {Theorem}  Suppose that $\varphi 
\vert BW_{n-1}\to MT_{n-1}$ is injective.\\
 Then $\varphi \vert BW_{2k}^{(0)}\to 
MT_{2k}^{(0)}$ is injective, when $n=2k$,\\
 and 
 $\varphi \vert BW_{2k+1}^{(1)}\to 
MT_{2k+1}^{(1)}$ is injective, when 
$n=2k+1$.\end{Theorem}   
\begin {Proof}  
In the case $n=2k$ we know that $\varphi \vert 
BW_{k}$ is an isomorphism to  $MT_{k}$. We may 
then choose elements $t_{c}\in BW_{k}, c\in 
C_{k}$, spanning $BW_{k}$, with $\varphi (t_{c})$ 
represented by a totally descending tangle 
$T_{c}$ say, having connector $c$.

By  corollary  \ref{5.6a} we have  $BW_{2k}^{(0)} 
=BW_{k}f_{k}BW_{k}$. This is spanned by $ \left 
\vert C_{k}\right \vert ^{2}$ elements 
$t_{c}f_{k}t_{d},\,  c,d\in C_{k}$. It is 
enough to prove that the images of these elements 
are independent in $MT_{2k}$.   

Now these images are represented by  the tangles 
$T_{c}F_{k}T_{d}$. Different pairs of connectors 
$(c,d)$ give tangles $T_{c}F_{k}T_{d}$ with 
different connectors in $C_{2k}$, since the 
tangles consist of a top and a bottom half, each 
with $k$ arcs, affected independently by the 
connectors $c$ and $d$. The tangles then 
represent independent elements in $MT_{2k}$, by 
 corollary  \ref{4.2a}.

Similarly when $n=2k+1$ we know that $\varphi 
\vert BW_{k+1}$ is an isomorphism to  
$MT_{k+1}$. We may then choose spanning elements 
$t_{c}\in BW_{k+1}, c\in C_{k}$,  with $\varphi 
(t_{c})$ represented by a totally descending 
tangle $T_{c}$ say, having connector $c$. Again, 
by  corollary  \ref{5.6a}, we have a spanning set 
$\lbrace t_{c}S(f_{k})t_{d}\rbrace ,\,  
c,d\in C_{k+1}$  with $ \left \vert C_{k+1}\right 
\vert ^{2}$ elements, for the ideal 
$BW_{2k+1}^{(1)}$. 

The images of these elements are represented by 
tangles $T_{c}S(F_{k})T_{d}$. Once more we can 
see that different pairs of connectors $(c,d)$ 
give tangles with different connectors in 
$C_{2k+1}$ because all but one of the arcs stays 
either in the top or in the bottom of the tangle. 
This guarantees independence in $MT_{2k+1}$, as 
before.\end{Proof} 

 \begin {Rem}  We could in fact show that the 
composite tangles used in this proof are 
themselves totally descending, for some suitable 
ordering of their arcs.\end{Rem}   
We continue in the next section to examine 
$BW_{n}^{(r)}$ for larger  $r$ having established 
here the start of our induction on $r$. Note that 
we could prove similarly that 
$BW_{2k+r}^{(r)}=BW_{k+r}S^{r}(f_{k})BW_{k+r}$ 
and find a spanning set of $ \left \vert 
C_{k+r}\right \vert ^{2}$ elements. However, a 
similar attempt to prove the (false) result for  
$r>1$ that these are independent would fail, 
because some different pairs of connectors in 
$C_{k+r}$ can yield the same connector in 
$C_{2k+r}$.

\section {Isomorphism between Kauffman's tangle algebras and the
Birman-Wenzl algebras}

We finish the proof of injectivity of $\varphi :BW_n\to MT_n$ 
by proving the remaining induction step, namely 
that if $\varphi \vert BW_{n-1}$ is injective, 
and $\varphi \vert BW_{n}^{(r-2)}$ is 
injective then $\varphi \vert BW_{n}^{(r)}$ is 
injective. We do this by finding a complementary 
subspace $V_{n}^{(r)}$ to $ BW_{n}^{(r-2)}$  in $ 
BW_{n}^{(r)}$  on which $\varphi $ is injective.

We recall the definition of $V_{n}^{(r)}$ given 
in section 4 as the subspace spanned by $\lbrace 
b_{\pi }BW_{2k}^{(0)}b_{\tau }b_{\mu }\rbrace $, 
where $n=2k+r$,  $\alpha (b_{\pi }),\, b_{\mu }$ are $(2k,r)$ Lorenz braids in 
$B_{n}$ and $b_{\tau }$ is a positive permutation 
braid on the last $r$ strings in $S^{2k}(B_{r})$. 
Following the scheme of proof in theorem \ref{4.5} we 
already know, by induction on $n$, that $\varphi 
\vert V_{n}^{(r)}$ is injective.

 It remains to 
show that $V_{n}^{(r)}+ BW_{n}^{(r-2)} = 
BW_{n}^{(r)}$.  Since $ BW_{n}^{(r)}$ is the 
2-sided ideal generated by $f_{k}$, and $f_{k}\in 
V_{n}^{(r)}$ we need only show that $ 
V_{n}^{(r)}+ BW_{n}^{(r-2)} $ is a 2-sided ideal. 
Now $\alpha (V_{n}^{(r)})=V_{n}^{(r)}$, since the 
elements $b_{\tau }$ in $S^{2k}(BW_{r})$ commute 
with $BW_{2k}$. Hence it is enough to show that $ 
V_{n}^{(r)}+ BW_{n}^{(r-2)} $ is a left ideal.
\begin {Proposition}\label{6.1}  Let $n=r+2k$ and let 
$X_{n}^{(r)}$ be the subspace spanned by the set $\lbrace 
b_{\pi }b_{\tau }BW_{2k}f_{k}\rbrace $, where 
$\alpha (b_{\pi })$ is a  $(2k,r)$ Lorenz braid 
and $b_{\tau }$ is a positive permutation braid 
in $S^{2k}(B_{r})$. Suppose also that $\varphi 
\vert BW_{n-1}$ is injective and that $r\geq 
2$. Then
$$
	L_{n}^{(r)}=X_{n}^{(r)}+BW_{n}^{(r-2)}
$$
is a left ideal.\end{Proposition} 
\begin {Corollary}   $ V_{n}^{(r)}+ 
BW_{n}^{(r-2)} $ is a left ideal, under the 
hypotheses of proposition \ref{6.1}, and hence theorem 
\ref{4.5} is established.\end{Corollary}   
\begin {Proof}  Since $L_{n}^{(r)}$ is a left 
ideal, by \ref{6.1}, it follows that $ V_{n}^{(r)}+ 
BW_{n}^{(r-2)} $ is a left ideal, by noting that 
$BW_{2k}^{(0)}=BW_{2k}f_{k}BW_{2k}$. \end{Proof}  
The proof of proposition \ref{6.1} 
occupies the rest of this section. The principal 
ingredient is an analysis of the elements 
$g_{i}b_{\pi }$ and $e_{i}b_{\pi }$ for positive 
permutation braids $b_{\pi }$. The following two 
lemmas are a consequence primarily of the braid 
relations.
\begin {Lemma} \label{6.2} Let $\rho $ be any 
permutation, and let $\rho _{1}$ be the 
permutation $\rho \circ (i\, i+1)$. Then 
the positive permutation braid $b_{\rho _1}$ 
satisfies the equation 
\begin{eqnarray*}  
 b_{\rho _1}=g_{i}b_{\rho }&\quad &\mbox{if } \rho (i)<\rho 
(i+1),\\   b_{\rho }=g_{i}b_{\rho_1} 
&\quad &\mbox{if } \rho
(i)>\rho (i+1).
\end{eqnarray*}
\end{Lemma}   
\begin {Proof}  If $\rho (i)<\rho (i+1)$ then 
each pair of strings in the braid $g_{i}b_{\rho 
}$ crosses at most once, so it is a positive 
permutation braid. Its permutation is $\rho 
_{1}$, so $g_{i}b_{\rho }=b_{\rho _1}$. 

If $\rho (i)>\rho (i+1)$ then $\rho _1(i)<\rho 
_1(i+1)$ and the same argument holds with $\rho 
_1$ in place of $\rho $. \end{Proof} 
\begin {Corollary}  Any positive permutation 
braid $b_{\rho }$ can be written as the product 
of a word in $\lbrace g_{i}\rbrace ,\, i\mathbin{\not =}\ell $, and an $(\ell ,r)$ 
Lorenz braid.\end{Corollary}   
\begin {Proof}  By induction on the length of 
$b_{\rho }$, using \ref{6.2} to write $b_{\rho 
}=g_{i}b_{\rho _1}$ for some $i\mathbin{\not 
=}\ell $ if $b_{\rho }$ is not already an $(\ell 
,r)$ Lorenz braid.\end{Proof} 
\begin {Lemma}\label{6.3}  Let $\rho $ be any 
permutation with $\rho (i+1)=\rho (i)+1$. Then 
$g_{i}b_{\rho }=b_{\rho }g_{\rho (i)}$ and 
$e_{i}b_{\rho }=b_{\rho }e_{\rho (i)}$. 
\end{Lemma}  
\begin {Proof}  This can be viewed as allowing us 
to pass a simple crossing along two parallel 
strings from top to bottom of a braid. By the 
hypothesis on $\rho $, both $g_{i}b_{\rho }$ and 
$b_{\rho }g_{\rho (i)}$ are positive permutation 
braids, and both have the same permutation. Hence 
they are equal, using only the braid relations, 
by the fundamental theorem on positive 
permutation braids. It follows that 
$g_{i}^{-1}b_{\rho }=b_{\rho }g_{\rho (i)}^{-1}$ 
and hence, by the skein relation, that $z\, e_{i}b_{\rho }=z\, b_{\rho }e_{\rho 
(i)}$. 

The lemma follows, if we assume that $z$ is 
invertible in $\Lambda $.  Without inverting $z$ 
the result follows by induction on the length of 
$b_{\rho }$, together with the relation 
$e_{i}g_{i+1}g_{i}=g_{i+1}g_{i}e_{i+1}$ and its 
reverse in $BW_{n}$. For we can write $b_{\rho 
}=g_{j}b_{\rho _1}$ for some $j$. Then 
$j\mathbin{\not =}i$, since the strings $i$ and 
$i+1$ do not cross under $\rho $. 

If $j=i+1$ then $\rho (i+2)<\rho (i+1)=\rho 
(i)+1$, so $\rho (i+2)<\rho (i)$. We can then, by 
lemma \ref{6.2}, write $b_{\rho }=g_{i+1}g_{i}b_{\rho 
_2}$, and then $e_{i}b_{\rho 
}=g_{i+1}g_{i}e_{i+1}b_{\rho _2}$. Now $\rho 
_2(i+2)=\rho (i+1)=\rho _2(i+1)+1$ and $b_{\rho 
_2}$ is shorter than $b_{\rho }$, so that we can 
use induction.

A similar argument can be used when $j=i-1$, 
while otherwise $ \left \vert i-j\right \vert 
>2$, and $e_{i}g_{j}=g_{j}e_{i}$, giving an 
immediate inductive proof.\end{Proof} 
\begin {Lemma}\label{6.4}   
$X_{n}^{(r)}S^{2k}(BW_{r})\subset L_{n}^{(r)}$. 
\end{Lemma}   
\begin {Proof}   $X_{n}^{(r)}e_{j}\subset 
BW_{n}^{(r-2)}\subset L_{n}^{(r)}$ for $j>2k$, 
since $f_{k}e_{j}\in BW_{n}^{(r-2)}$ for $j>2k$. 

Let $b_{\tau }$ be any positive permutation braid 
in $S^{2k}B_{r}$ and let $j>2k$. Then by \ref{6.2}, 
either \begin{eqnarray*}b_{\tau }g_{j}&=&b_{\tau '} 
\\  \mbox {or } b_{\tau }g_{j}&=&b_{\tau 
'}g_{j}^{2}=b_{\tau '}+zb_{\tau }-zb_{\tau 
}e_{j}. \end{eqnarray*}
Hence $xg_{j}\in L_{n}^{(r)}$ for any spanning 
element $x=b_{\pi }w_{2k}f_{k}b_{\tau }\in 
X_{n}^{(r)}$. 

Thus $X_{n}^{(r)}g_{j}\subset L_{n}^{(r)}$ for 
$j>2k$. \end{Proof} 

We now continue the proof of \ref{6.1}, to show that  
$L_{n}^{(r)}$ is a left ideal. Lemma \ref{6.4} shows in 
particular that $L_{n}^{(r)}b_{\tau }\subset 
L_{n}^{(r)}$ for $b_{\tau }\in S^{2k}(B_{r})$. It 
is enough to show that $e_{i}x,\, g_{i}x\in L_{n}^{(r)}$ for each $x=b_{\pi 
}w_{2k}f_{k}\in X_{n}^{(r)}$ and each $i$, where 
$\alpha (b_{\pi })$ is a  $(2k,r)$ Lorenz braid 
and $w_{2k}\in BW_{2k}$. 

Suppose then that $x$ and $i$ are given. We may 
further suppose that $\pi (i+1)>\pi (i)$, 
otherwise $b_{\pi }=g_{i}b_{\pi _1}$ with $\pi 
_1(i+1)>\pi _1(i)$. We then need only prove that 
$e_{i}x'\in L_{n}^{(r)}$ where $x'=b_{\pi 
_1}w_{2k}f_{k}$,  since 
$g_{i}x=g_{i}^{2}x'=x'+zg
_{i}x'-zg_{i}e_{i}x'=x'+zx-\lambda ze_{i}x'$ and 
$e_{i}x=e_{i}g_{i}x'=\lambda e_{i}x'$ from the 
skein and delooping relations. 

Since $\pi $ is a reverse $(2k,r)$ Lorenz 
permutation then $\pi (i+1)=\pi (i)+1$ if either 
$\pi (i+1)\leq 2k$ or $\pi (i)>2k$.   By \ref{6.3} 
$e_{i}x=b_{\pi }e_{\pi (i)}w_{2k}f_{k}$ in either 
case. This lies in $X_{n}^{(r)}$ if $\pi (i)<2k$ 
and in $BW_{n}^{(r-2)}$ if $\pi (i)>2k$, and 
similarly $g_{i}x\in L_{n}^{(r)}$. It remains to 
deal with $e_{i}x$ and $g_{i}x$ when $\pi 
(i+1)>2k$ and $\pi (i)\leq 2k$. In this case 
$g_{i}b_{\pi }$ is a reverse $(2k,r)$ Lorenz 
braid, by \ref{6.2}, so that $g_{i}x\in X_{n}^{(r)}$ 
and we are left to consider $e_{i}x$. 

Given $\pi $ and $i$, let $\rho $ be the 
permutation given by
$$\rho (j)=\cases{2k& $j=\pi (i),$   \cr 
 j-1,& $\pi (i)<j\leq 2k,$   \cr  j+1,& 
$2k+1\leq j<\pi (i+1),$   \cr  2k+1,& 
$j=\pi (i+1),$   \cr  j,&otherwise. 
 \cr  }$$
Now  $\rho $ only makes pairs of strings cross 
which have not already been made to cross by the 
reverse Lorenz braid $b_{\pi }$, so that $b_{\pi 
}b_{\rho }$ is also a positive permutation braid. 
Then $b_{\pi }b_{\rho }=b_{\pi _1}$, where $\pi 
_1=\rho \circ \pi $. Now $\rho $ permutes the 
first $2k$ strings and the last $r$ strings among 
themselves, moving $\pi (i)$ to $2k$ and $\pi 
(i+1)$ to $2k+1$,  so  $x=b_{\pi _1}(b_{\rho 
})^{-1}w_{2k}f_{k}$ with $b_{\rho }\in 
BW_{2k}S^{2k}(B_{r})$. Note that $\pi 
_{1}^{-1}(i)<\pi _{1}^{-1}(j)$ for $i<j\leq 
2k-1$. 

It is enough, by lemma \ref{6.4}, to show that 
$e_{i}x'\in L_{n}^{(r)}$, for $x'=b_{\pi 
_1}w'_{2k}f_{k}$. Now $\pi _1(i+1)=2k+1=\pi 
_1(i)+1$ so, by \ref{6.2},  $e_{i}x'=b_{\pi 
_1}e_{2k}w_{2k}'f_{k}$. This does not finish the 
proof, since the element $e_{2k}$ is stuck 
between $BW_{2k}$ and $S^{2k}(BW_{r})$ and we 
have to use our inductive knowledge of 
$w'_{2k}f_{k}\in BW_{2k}^{(0)}$ to free it.
\begin {Lemma}\label{6.6}  Suppose that $\varphi 
\vert BW_{2k}^{(0)}$ is injective. Then every 
element in $BW_{2k}f_{k}$ is a linear combination 
of elements in  the sets $$g_{m}g_{m+1}\dots{} 
g_{2k-2}e_{2k-1}BW_{2k}f_{k},  m=1,\dots{} 
,{2k-2}, \mbox {and } e_{2k-1}BW_{2k}f_{k}.$$ \end{Lemma}   
\begin {Lemma}\label{6.7}  For each $ m=1,\dots{} 
,{2k-2} $ and each positive permutation braid 
$b_{\rho }$ with $\rho ^{-1}(i)<\rho ^{-1}(j)$ 
for $i<j\leq 2k-1$ and $\rho ^{-1}(2k+1)=\rho 
^{-1}(2k)+1$  we have
$$b_{\rho }e_{2k} g_{m}g_{m+1}\dots{} 
g_{2k-2}e_{2k-1}=b_{\rho '}e_{2k-1}$$
for some positive permutation braid $b_{\rho '}$. 
\end{Lemma}  

\medskip 
Proposition \ref{6.1} then follows from \ref{6.6} and \ref{6.7}, 
since we can write the element $e_{i}x'$ as a 
linear combination of elements of the form 
$b_{\rho '}BW_{2k}f_{k}$. All of these lie in 
$L_{n}^{(r)}$, since any positive permutation 
braid $b_{\rho '}$ can be written as the product 
$b_{\pi }b_{\rho ''}$ of a reverse $(2k,r)$ 
Lorenz braid  $b_{\pi }$ with a positive braid 
which does not involve the generator $g_{2k}$, by 
the corollary to lemma \ref{6.2}, applied to the 
reverse braids.

\medskip \begin {Proof}  [of lemma \ref{6.6}]  
By hypothesis, $\varphi $ gives an isomorphism 
from $BW_{2k}f_{k}\subset BW_{2k}^{(0)}$ to 
$MT_{2k}F_{k}$. Now every element of 
$MT_{2k}F_{k}$ can be written as a linear 
combination of totally descending tangles 
$T_{c}$, where the connectors $c$ join points of 
the top to the top in some way, and join the 
bottom points as for $F_{k}$.   We may choose the 
order of strings for each connector $c$ as we 
wish, so let us assume that in each tangle 
$T_{c}$ the string whose end point is at position 
$2k$ on the top lies above all the others. By 
isotopy of the strings we may then write each of 
these tangles $T_{c}$ as
$$G_{m}G_{m+1}\dots{} G_{2k-2}E_{2k-1}TF_{k}$$ 
for some $m=1,\dots{} ,2k-2$ and some $ T\in 
MT_{2k} $ as illustrated in figure 6.1. The 
isomorphism $\varphi $ then gives a spanning set 
for $BW_{2k}f_{k}$ as stated.\end{Proof} 

\begin{center}
\BWsixa\\
{Figure 6.1}
\end{center}

 \begin {Proof} [ of lemma \ref{6.7}]  We 
have $$ b_{\rho }e_{2k} g_{m}g_{m+1}\dots{} 
g_{2k-2}e_{2k-1}= b_{\rho }g_{m}g_{m+1}\dots{} 
g_{2k-2}e_{2k}e_{2k-1}. $$ Now the reverse braid 
$g_{2k-2}\dots{} g_{m+1}g_{m}\alpha (b_{\rho })$ 
is a positive permutation braid,  $b_{\rho _1}$ 
say, since $g_{2k-2}\dots{} g_{m}$ is a positive 
permutation braid on the first $2k-1$ strings 
only, while $\alpha (b_{\rho })=b_{\rho ^{-1}}$ 
does not make these strings cross. Now $\rho 
_{1}(2k+1)=\rho _{1}(2k)+1$, so either 
$g_{2k}g_{2k-1}b_{\rho _1}$ or 
$g_{2k}^{-1}g_{2k-1}^{-1}b_{\rho _1}$ is a 
positive permutation braid, $b_{\rho _2}$, say, 
depending on whether $\rho _{1}(2k-1)<\rho 
_{1}(2k)$ or $\rho _{1}(2k-1)>\rho _{1}(2k)$, by 
\ref{6.2}.

We can write $e_{2k-1}e_{
2k}=e_{2k-1}g_{2k}g_{2k-1
}=e_{2k-1}g_{2k}^{-1}g_{2k-1}^{-1}$ by the 
relations in  $BW_{n}$. Then 
$e_{2k-1}e_{2k}b_{\rho _1}=e_{2k-1}b_{\rho _2}$. 
Apply the reversing map to give 
\begin{eqnarray*} 
b_{\rho }e_{2k} g_{m}g_{m+1}\dots{} 
g_{2k-2}e_{2k-1}&=&\alpha ( e_{2k-1}e_{2k}b_{\rho 
_1}) \\  &=&\alpha (e_{2k-1}b_{\rho 
_2}) \\  &=&b_{\rho ''}e_{2k-1}, 
\end{eqnarray*}
 where $\rho ''=\rho 
_{2}^{-1}$. \end{Proof} 

\medskip 
This concludes the proof of proposition \ref{6.1}, and 
the inductive proof of theorem \ref{4.5}. We have now 
established that $\varphi $ is an isomorphism 
from $BW_{n}$ to $MT_{n}$ for all $n$, so that we 
are able to use tangle based arguments in dealing 
with the algebra $BW_{n}$. We have established 
its dimension over $\Lambda $ and also the 
geometric description of the natural chain of 
ideals generated by the elements $f_{k}$, so we 
can also study the composition series of this 
chain by using the corresponding ideals in 
$MT_{n}$ generated by $F_{k}$.

\bigskip

%

\end{document}